\documentclass[a4,reqno]{amsart}

\usepackage{graphicx}
\usepackage{pgf,tikz,float}
\usetikzlibrary{arrows}
\usepackage{hyperref}

\hypersetup{colorlinks, citecolor=[rgb]{0.0,0.3,0.1}, linkcolor=[rgb]{0.3,0.0,0.1},
	urlcolor=[rgb]{0.0,0.0,0.4}}

\usepackage[nobysame]{amsrefs}

\allowdisplaybreaks

\author{Ferenc Fodor}

\title[Perimeter approximation of convex discs]{Perimeter approximation of convex discs in the hyperbolic plane and on the sphere}

\address{Department of Geometry, Bolyai Institute, University of Szeged, Aradi v\'ertan\'uk tere 1, 6720 Szeged, Hungary}

\email{fodorf@math.u-szeged.hu}



\newcommand{\R}{\mathbb{R}}
\DeclareMathOperator{\per}{Per}
\DeclareMathOperator{\dev}{dev}
\DeclareMathOperator{\Pe}{{\mathcal P}}
\DeclareMathOperator{\bd}{{\mathrm bd}}
\DeclareMathOperator{\HH}{\mathbb{H}}
\DeclareMathOperator{\Sp}{\mathbb{S}}

\newtheorem{theorem}{Theorem}[section]

\newtheorem{corollary}[theorem]{Corollary}

\begin{document}

\begin{abstract}
Eggleston \cite{Eggleston-1957} proved that in the Euclidean plane the best approximating convex $n$-gon to a convex disc $K$ is always inscribed in $K$ if we measure the distance by perimeter deviation. We prove that the analogue of Eggleston's statement holds in the hyperbolic plane, and we give an example showing that it fails on the sphere. 
\end{abstract}

\maketitle

\section{Introduction and main results}
We call a compact, convex set $K\subset\R^2$ whose interior is non-empty a {\em convex disc}.
The perimeter of $K$ is denoted by $\per (K)$. Let $K$ and $L$ be both convex discs. 
The {\em perimeter deviation of $K$ and $L$} is defined as
$$\dev (K, L)=\per(K\cup L)-\per(K\cap L).$$
We note that although the perimeter deviation is often used to measure the distance of convex figures, it does not define a proper metric on the set of all
convex discs. For another notion of perimeter deviation, which is in fact a metric, see, for example, Florian \cite{Florian} and the references therein.

Eggleston \cite{Eggleston-1957}, among other questions, investigated 
how well a convex disc can be approximated
by convex polygons of a given number of vertices in the sense of perimeter deviation. 
For a positive integer $n\geq 3$, let $\Pe(n)$ denote the set of convex polygons with at most
$n$ vertices. Let 
$$\delta_{\dev}(K,n)=\inf\{\dev(K, P): P\in\Pe(n)\}.$$
A simple compactness argument shows that for each convex disc $K$ and positive
integer $n\geq 3$, there exists a $P\in\Pe(n)$ which minimizes the 
perimeter deviation from $K$, that is, $\dev(K,P)=\delta_{\dev}(K,n)$. 

Eggleston proved the following beautiful statement, cf. \cite[Lemma~4 on p. 353]{Eggleston-1957}.
\begin{theorem}[Eggleston, 1957]\label{Eggleston}
Let $K$ be a convex disc and $n\geq 3$ a positive integer. If $P\in\Pe(n)$ is
such that $\dev(K,P)=\delta_{\dev}(K,n)$, then $P$ is inscribed in $K$, that is,
$P\subset K$ and the vertices of $P$ are on the boundary of $K$.  
\end{theorem}

According to a classical result of Dowker \cite{Dowker-1944}, 
the minimum area of convex $n$-gons containing a given convex disc $K$ 
is a convex function of $n$, 
and the maximum area of convex $n$-gons contained in $K$ 
is a concave function of $n$. This result was later extended 
for perimeter in place of area by L. Fejes T\'oth \cite{LFT-1955}, Eggleston \cite{Eggleston-1957},
and Moln\'ar \cite{Molnar-1955}, independently from each other. 
Thus, it follows from Theorem~\ref{Eggleston} that for a fixed convex disc $K$, 
the minimum perimeter deviation of convex $n$-gons from $K$ is also
a concave function of $n$.  

Let $\HH^2$ denote the hyperbolic plane, and for two points $p,q\in\HH^2$ let $pq$ denote the (hyperbolic) segment with endpoints $p$ and $q$. The hyperbolic distance $d_H(p,q)$ of $p$ and $q$ is
the length of $pq$.  

Let $\Pe_H(n)$ denote the set of all
convex polygons in $\HH^2$ with at most $n$ vertices for $n\geq 3$. 
Similarly to the Euclidean case, we define
$$\delta_{\dev_H}(K,n)=\inf\{\dev (K,P):P\in\Pe_H(n)\}.$$ 
For any fixed $K$ and positive integer $n\geq 3$, there exits a convex 
polygon $P\in\Pe_H(n)$ such that $\dev_H(K,P)=\delta_{\dev_H}(K,n)$.
In this paper we extend Theorem~\ref{Eggleston} to the hyperbolic plane $\HH^2$ as follows.
\begin{theorem}\label{main}
Let $K$ be a convex disc in $\HH^2$ and $n\geq 3$ a positive integer. 
If $P\in\Pe_H(n)$ is such that $\dev(K,P)=\delta_{\dev_H}(K,n)$, 
then $P$ is inscribed in $K$, that is,
$P\subset K$ and the vertices of $P$ are on the boundary of $K$.  
\end{theorem}

The analogues of Dowker's theorem both for area and perimeter also hold on the sphere  $\Sp^2$ and
the hyperbolic plane $\HH^2$. These were proved by Moln\'ar \cite{Molnar-1955} and L. Fejes T\'oth \cite{LFT-1958}.
Thus, Theorem~\ref{main}, combined with the hyperbolic version of Dowker's theorem for the 
maximum perimeter of convex (hyperbolic) $n$-gons contained in a given convex disc $K$, 
implies the following statement.

\begin{corollary}
The minimum perimeter deviation of convex $n$-gons from a given 
convex disc $K$ is a concave function of $n$ in the hyperbolic plane $\HH^2$.  
\end{corollary}

On the unit sphere $\mathbb{S}^2$, the distance of two non-antipodal points $p,q$ is the length of the shorter arc of the unique great circle through $p$ and $q$. The distance of two antipodal points is $\pi $. We call a closed set $K$ on $\mathbb{S}^2$ (spherically) convex if it is contained in an open hemisphere and for any $p,q\in K$, the shorter arc of the unique great circle connecting $p$ and $q$ is also contained in $K$.   
One may naturally define the perimeter deviation $\dev_S(K,L)$ of two  convex discs $K, L$ on the unit sphere as in the Euclidean plane and hyperbolic plane. Again, for a convex disc $K$ and $n\geq 3$, there exists a convex spherical polygon $P$ with at most $n$ vertices such that $\dev_S(K,P)=\delta_{\dev_S}(K, n)$. However, $P$ may not necessarily be contained in $K$ (or contain $K$) as shown by an example in Section~\ref{sphere}. 

\section{Proof of Theorem~\ref{main}}
In this section we work in the hyperbolic plane $\HH^2$, thus all notions, such as distance, convexity, perimeter, perimeter deviation, etc. are always understood in the hyperbolic sense without mentioning this fact explicitly. We think of $\HH^2$ as a $2$-dimensional Riemannian manifold of constant curvature $-1$, such as the Beltrami-Klein model, see more on this below. By the {\em curvature} of a $C^2$ curve in $\HH^2$ we mean its geodesic curvature. A compact set $K\subset\HH^2$ is (geodesically) convex, if for any $x,y\in K$, the geodesic segment $xy$ is contained in $K$. We note that there are other forms of convexity in $\HH^2$, for example, $h$-convexity where we require that the whole region bounded by the two horocyclic arcs connecting $x$ and $y$ is contained in $K$, or $\lambda$-convexity where one requires that $K$ contains the region bounded by two congruent hypercyclic arcs of radius $\lambda$ through $x$ and $y$. For more information we refer to \cite{GR99}.

We essentially follow a similar but somewhat more complicated argument to that of Eggleston in \cite{Eggleston-1957}.

First, note that the set of all compact, convex sets forms a complete metric space with respect to the 
Hausdorff distance in $\HH^2$. Furthermore, the perimeter deviation function is continuous 
on this space. Thus, it is enough to prove the theorem on a suitable dense subspace of convex discs. 
We select this dense subspace the following way: We assume that the boundary $\bd K$ of $K$ is $C^2_+$ smooth, meaning that it is twice continuously differentiable at every point and the geodesic curvature is strictly positive everywhere. 

We start the proof by examining the difference between the length of a chord and the corresponding arc of $\bd K$ cut off by the chord. 

In the following argument, we will work in the Beltrami-Klein model $D$ of the hyperbolic plane $\mathbb{H}^2$ whose points are the interior points of the unit radius circular disc $B^2$ centred at the origin, and whose lines (geodesics) are the Euclidean open line segments with endpoints on the boundary $S^1$ of $B^2$. We define the distance of two points $p, q\in D$ as
$$d_H(p,q)=\frac 12|\ln (abpq)|,$$
where $a$ and $b$ are the intersection points of the line $pq$ with $S^1$ such that $a$ is on the side of $p$ and $b$ is on the side of $q$. The symbol $(abpq)$ denotes the cross-ratio of the points $a,b,p,q$ in this order. It is well-known that the Gaussian curvature of this model is constant $-1$, with this particular metric. Now, if $p(x,y)$ is a point of $D$, where $x$ and $y$ are its Euclidean coordinates in a Cartesian coordinate system centred at the origin, then the hyperbolic coordinates of $p(x_h, y_h)$ are the following
$$x_h=\frac 12\ln\frac{1+x}{1-x}, \quad y_h=\frac 12\ln \frac{\sqrt{1-x^2}+y}{\sqrt{1-x^2}-y}.
$$ 
The first fundamental form of $D$ is 
$$ds^2=\frac{(1-y^2)dx^2+2xydxdy+(1-x^2)dy^2}{(1-x^2-y^2)^2},$$
see, for example, \cite{Cannon}.

Let $K\subset D$ be a (geodesically) convex disc in $D$ whose boundary is $C^2_+$ smooth. Since geodesic segments in $D$ are exactly the Euclidean segments, the disc $K$ is convex in the hyperbolic sense exactly if it is convex in the Euclidean sense. Assume that $o\in \bd K$ and that the $x$-axis supports $K$ at $o$. Then, in a suitably small neighbourhood of $o$, the boundary of $K$ can be represented by a convex function $f$ such that $f(x)=(\kappa/2) x^2+o(x^2)$ as $x\to 0$, and $\kappa>0$. A standard calculation shows that the geodesic curvature 
of $\bd K$ at $o$ is $\kappa$. 

For sufficiently small $x$, let $s(x)$ denote the arc-length of $\bd K$ between $o$ and the point $(x,f(x))$. Then
\begin{align}
s(x)&=\int_0^{x}\frac{((1-f^2(\tau))+2\tau f(\tau)f'(\tau)+(1-\tau^2)(f'(\tau)^2))^{1/2}}{1-\tau^2-f^2(\tau)}d\tau\notag\\
&=\int_0^{x}\frac{\sqrt{1+\kappa^2\tau^2+o(\tau^2)}}{1-\tau^2+o(\tau^2)}d\tau.\label{elso}
\end{align}
After substituting the Taylor series of $\sqrt{1+z}$ around $z=0$ and that of $(1-z)^{-1}$ around $z=0$ in \eqref{elso}, we obtain

\begin{align}
s(x)&=\int_0^{x}\left (1+\frac{\kappa^2\tau^2}{2}+o(\tau^2)\right )\left (1+\tau^2+o(\tau^2)\right )d\tau\notag\\
&=\int_0^{x}1+\frac{\kappa^2+2}{2}\tau^2+o(\tau^2)d\tau\notag\\
&=\left (\tau+\frac{\kappa^2+2}{6}\tau^3+o(\tau^3)\right ]_0^{x}\notag\\
&=x+\frac{\kappa^2+2}{6}x^3+o(x^3)\text{ as }x\to 0^+.\label{masodik}
\end{align}

First, let $l=l(\delta)$ be the line with Euclidean equation $y=\delta$. For sufficiently small $\delta>0$, the line $l$ intersects $\bd K$ at $x_+(\delta)>0$ ($x_-(\delta)<0$) such that $f(x_+(\delta))=\delta$ ($f(x_-(\delta))=\delta$). Due to the definition of $f$, $x_+(\delta)=\sqrt{2/\kappa}\delta^{1/2}+o(\delta^{1/2})$
($x_-(\delta)=-\sqrt{2/\kappa}\delta^{1/2}+o(\delta^{1/2})$) as $\delta\to 0^+$.

Thus, by \eqref{masodik}, the arc of $\bd K$ between $o$ and the positive intersection point of $l$ and $\bd K$ has length
\begin{equation}\label{negyedik}
s(x_+(\delta))=x_+(\delta)+\frac{\kappa^2+2}{6}x^3_+(\delta)+o(x^3_+(\delta))\text{ as }\delta\to 0^+.
\end{equation}
Clearly, a similar formula holds for the length of the arc of $\bd K$ between the intersection point with (negative) $x$-coordinate $x_-(\delta)$ and $o$. 

The (hyperbolic) length of the segment between the $y$-axis and the (positive) intersection point with $\bd K$ is the following
\begin{align}
s_l(\delta)&=\int_0^{x_+(\delta)}\frac{\sqrt{1-\delta^2}}{1-x^2-\delta^2}dx\notag\\
&=\frac 12\ln \frac{\sqrt{1-\delta^2}+x_+(\delta)}{\sqrt{1-\delta^2}-x_+(\delta)}\notag\\
&=x_+(\delta)+\frac 13x^3_+(\delta)+o(\delta^2x_+(\delta))\text{ as }x_+(\delta)\to 0^+\label{harmadik},
\end{align}
and, again, a similar formula holds for the length of the segment between the negative intersection point of $l$ and $\bd K$ and the $y$-axis.

From \eqref{negyedik} and \eqref{harmadik}, and the expressions of $x_+(\delta)$ and $x_-(\delta)$, we obtain that the difference of the arc of $\bd K$ and the chord at height $\delta$ is 
\begin{equation}\label{par-chord}
\frac{\kappa^2}{3}(x^3_+(\delta)+x^3_-(\delta))+o(x^3_+(\delta))+o(x^3_-(\delta))=O(\delta^{3/2}) \text{ as }\delta\to 0^+.
\end{equation}

Second, we assume that the Euclidean equation of the line $l$ is $y=\tan\theta\cdot x$, meaning that $l$ passes through $o$ and makes an angle $\theta$ with the positive part of the $x$-axis. If $\theta>0$ is sufficiently small, then for the $x$-coordinate $x(\theta)$ of the intersection point of $l$ and $\bd K$, different from $o$, the following holds
$$f(x(\theta))=\tan\theta\cdot x(\theta),$$ 
from which we obtain that
$$x(\theta)=2\tan\theta/\kappa+o(\tan\theta)=2\theta/\kappa+o(\theta)\text{ as } \theta\to 0^+.$$

Substituting $x(\theta)$ in \eqref{masodik}, we get that the arc-length of $\bd K$ between $o$ and the other intersection point of $l$ and $\bd K$ is
\begin{align}
s(\theta)=x(\theta)+\frac{\kappa^2+2}{6}x^3(\theta)+o(x^3(\theta))\label{ferde-egy}\text{ as }x(\theta)\to 0^+.
\end{align}	

At the same time, the length of the segment $l\cap K$ is 
\begin{align}
s_l(\theta)&=\tanh^{-1}(\sqrt{x^2(\theta)+f^2(x(\theta))})\notag\\
&=\tanh^{-1}(\sqrt{x^2(\theta)+\tan^2\theta x^2(\theta)})\notag\\
&=\tanh^{-1}(x(\theta)\sec \theta)\notag\\
&=x(\theta)\sec\theta+\frac 13x^3(\theta)\sec^3\theta+O(x^3(\theta)\sec^3\theta)\notag\\
&=x(\theta)+\frac 12x(\theta)\theta^2++\frac 13x^3(\theta)\sec^3\theta+O(x^3(\theta)\sec^3\theta).\label{ferde-ketto}
\end{align}
Now, by \eqref{ferde-egy} and \eqref{ferde-ketto}, the difference between the chord of $l$ and the corresponding part of $\bd K$ is
\begin{align}
s(\theta)-s_l(\theta)=\left (\frac{8\kappa^2}{6\kappa^3}-\frac{1}{\kappa}\right )\theta^2+o(\theta^3)=\frac {1}{3\kappa} \theta^2 +o(\theta^3)=O(\theta^2)\text{ as }\theta\to 0^+.\label{ferde-hossz}
\end{align}

The observations \eqref{par-chord} and \eqref{ferde-hossz} are elementary and probably well-known. We only included their detailed proofs because we could not find an explicit argument in the literature. 

Now, we turn to the actual proof of Theorem~\ref{main}.
Let $P\in\Pe_H(n)$ be an $n$-gon which minimizes the perimeter deviation from $K$, that is, 
$\dev(K,P)=\delta_{\dev_H}(K,n)$. We will denote the vertices of $P$ by $x_1,\ldots, x_n$ in 
a counter-clockwise cyclic order along $P$. It is clear that each side $x_ix_{i+1}$ 
has a common point with $K$, otherwise we could move 
it inwards and decrease the perimeter deviation using the monotonicity of perimeter in the hyperbolic plane, cf. \cite[Proposition~1.3]{Klain}. 
The assumption that $\bd K$ is $C^2_+$ yields that $K$ is strictly convex, that is, $\bd K$ contains no geodesic segment, and that $\bd K$ has a unique supporting line at each point, and therefore it cannot have vertices. 

The proof of Theorem~\ref{main} is indirect: we assume, on the contrary, that $P$ is not inscribed in $K$ and seek a contradiction. It is clear that if $P\subset K$, then the vertices of $P$ must be on $\bd K$, similarly to the Euclidean case, or otherwise we could increase the perimeter of $P$ by moving the vertices out to the boundary of $K$. Therefore, the indirect assumption yields that $P$ has a side with at least one endpoint outside of $K$.
There are several possibilities how this may happen. We treat each such case and show that they all contradict to the best approximation property of $P$. 

We use the following notation, similar to \cite[Section~2]{Eggleston-1957}. Let the vertex $x_i$ be outside of $K$. We denote the internal angle of $P$ at $x_i$ by $\alpha_i$. Let $b_i$ be the last common point of the side $x_{i-1}x_i$ and $\bd K$, and let $c_i$ be the first common point of $x_ix_{i+1}$ and $\bd K$ in the counter-clockwise direction along $P$. Let the angle of the tangent of $\bd K$ at $b_i$ and $x_{i-1}x_i$ be denoted by $\beta_i$, and the angle of the tangent of $\bd K$ at $c_i$ and $x_ix_{i+1}$ be $\gamma_i$. Then, clearly, $\alpha_i+\beta_i+\gamma_i<\pi$.

Let us first consider the case when $P$ has a side, say $x_1x_2$, such that both $x_1$ and $x_2$ are outside of $K$.  Let $\delta>0$ be small and let $h =h(\delta)$ be the hypercycle that is the equidistant curve from the line $x_1x_2$ at distance $\delta$ in the half-plane of $x_1x_2$ containing $P$. Let $x_1'$ and $x_2'$ be the intersection points of the sides $x_nx_1$ and $x_2x_3$ with $h(\delta)$, respectively. Assume that $\delta$ is so small that both $x_1'$ and $x_2'$ are outside of $K$. Then the $n$-gon $P'$ with vertices $x_1',x_2',x_3\ldots, x_n$, in this order, is contained in $P$, and the intersection of the side $x_1'x_2'$ and $K$ is of positive length. Let the feet of the perpendiculars from $x_1'$ and $x_2'$ to $x_1x_2$ be $x_1''$ and $x_2''$, respectively. 
Then $x_1''x_2''x_2'x_1'$ is a Saccheri quadrilateral. It is known that the line through the midpoints of the segments $x_1'x_2'$, and $x_1''x_2''$ is perpendicular to both lines and thus it cuts $x_1''x_2''x_2'x_1'$ into two congruent Lambert quadrilaterals. Using known trigonometric relations for Lambert quadrilaterals, we obtain that 
$$\sinh (d(x_1'',x_2'')/2)=\sinh (d(x_1',x_2')/2)\cosh \delta,$$ 
from which it follows that
$$d(x_1'x_2')=d(x_1'',x_2'')+O(\delta^2) \text{ as }\delta\to 0^+.$$
By hyperbolic trigonometry, we obtain for $i=1,2$ that 
$$\sinh d(x_ix_i')=\sinh \delta \csc\alpha_i,$$
thus
$$d(x_ix_i')=\delta \csc\alpha_i+O(\delta^3)\text{ as }\delta\to 0^+,$$
and
$$\sinh d(x_ix_i'')=-\tanh\delta \cot \alpha_i,$$
thus
$$d(x_ix_i'')=-\delta \cot\alpha_i+O(\delta^3)\text{ as }\delta\to 0^+.$$

If $x_1x_2$ is tangent to $K$ at a relative interior point $x'\in x_1x_2$, then 
let $a$ and $b$ denote the intersection points of the segment $x_1'x_2'$ with $\bd K$ so that $a$ is closer to $x_1'$. 
Then $d(x', ab)\leq\delta$. 
By the positivity of the geodesic curvature of $\bd K$ at $x'$ and by \eqref{par-chord}, it holds that the difference of the arc-length of $\bd K$ between $a$ and $b$ and the length of the segment $ab$ is $O(\delta^{3/2})$ as $\delta\to0$, and thus, by using the estimates obtained above, we get that
$$\dev (K, P')=\dev (K, P)-\delta (\csc\alpha_1+\csc\alpha_2+\cot\alpha_1+\cot\alpha_2)+
O(\delta^{3/2}) \text{ as }\delta\to 0^+.$$
Since $\csc\alpha+\cot\alpha\geq 0$ for any $\alpha\in (0,\pi)$, the coefficient of $\delta$ is negative in the above expression. This contradicts the minimality of $P$, and thus $P$ cannot have such a side. 

\begin{figure}
\begin{center}
\includegraphics[scale=0.35]{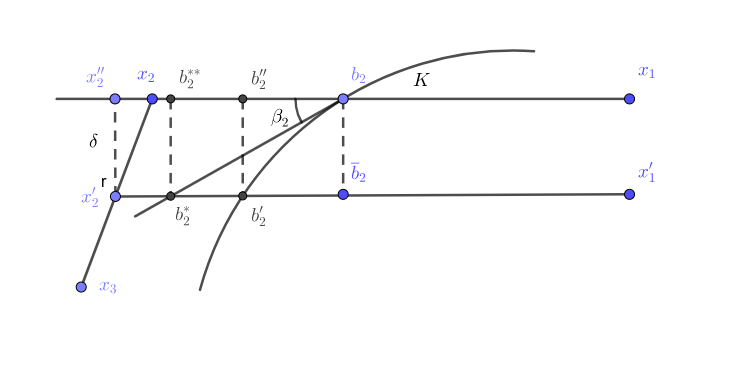}
\end{center}
	\caption{\label{fig:one}}
\end{figure}

If the side $x_1x_2$ cuts the boundary in two distinct points that are relatively interior to $x_1x_2$, then according to the previously introduced notation these intersection points are $c_1$ and $b_2$, and the tangents to $\bd K$ make an angle $\gamma_1$ and $\beta_2$ with $x_1x_2$, respectively. We introduce the following notations, see Figure~\ref{fig:one}. 
Let the last intersection point of $x_1'x_2'$ and $\bd K$ be $b_2'$. Let $b_2''$ be the perpendicular projection of $b_2'$ onto $x_1x_2$. Let $\overline{b}_2$ be the point on $x_1'x_2'$ whose perpendicular projection onto $x_1x_2$ is $b_2$. Let $b_2^*$ be the intersection point of the tangent line of $\bd K$ through $b_2$ and $x_1'x_2'$. Finally, let $b_2^{**}$ be the perpendicular projection of $b_2^*$ onto $x_1x_2$. 

We first note that, using hyperbolic trigonometry, we can conclude that 
$$d(b_2,b_2'')=d(\overline{b}_2,b_2')+O(\delta^2)\text{ as }\delta\to 0^+,$$
and similarly,
$$d(b_2'',b_2^{**})=d(b'_2,b_2^*)+O(\delta^2)\text{ as }\delta\to 0^+.$$
Now, let $\delta'=d(b_2^*,b_2^{**})$. Clearly, $\delta'<\delta$.
Similar as above, we obtain by hyperbolic trigonometry applied to the triangle $b_2b_2^*b_2^{**}$ that
$$\sinh d(b_2,b_2^*)=\sinh \delta'\csc \beta_2,$$
from where
$$d(b_2,b_2^*)=\delta'\csc\beta_2+O(\delta'^{3})\text{ as }\delta\to 0^+,$$
moreover,
$$\sinh d(b_2,b_2^{**})=-\tanh\delta'\cot\beta_2,$$
and 
$$d(b_2,b_2^{**})=-\delta'\cot\beta_2+O(\delta'^{3})\text{ as }\delta\to 0^+.$$

From trigonometric formulas for the corresponding Lambert quadrilateral we get that
$$\delta'=\delta+O(\delta^3)\text{ as }\delta\to 0^+.$$

Also, it is clear from the $C^2_+$ property of $\bd K$ that 
$$d(b_2', b_2^*)=O(d(b_2, b_2^*))=O(\delta^2)\text{ as }\delta\to 0^+,$$
and thus from all of the above,
$$d(b_2,b_2'')=d(b_2,b_2^{**})+O(\delta^2) \text{ as }\delta\to 0^+,$$
and
$$d(\overline{b}_2,b_2')=d(\overline{b}_2,b_2^*)+O(\delta^2) \text{ as }\delta\to 0^+.$$

Let $l(\delta)$ denote the length of the arc of $\bd K$ between $b_2$ and $b_2'$. From \eqref{masodik}, we obtain that 
$$l(\delta)-d(b_2,b_2^*)=O(\delta^3)\text{ as }\delta\to 0^+.$$

Finally, putting everything together, we obtain (similar to (19) in \cite{Eggleston-1957}) that
\begin{align*}
	\dev_H(K,P')=&\dev_H(K,P)-\delta (2\cot \beta_2-2\csc\beta_2+\csc\alpha_2+\cot\alpha_2\\
	&+2\cot \gamma_1-2\csc\gamma_1+\csc\alpha_1+\cot\alpha_1 )+O(\delta^2) \text{ as }\delta\to 0^+,
\end{align*}
and thus, by the optimality of $P$, it must hold that
\begin{equation}\label{second}
\cot\frac 12 \alpha_1+\cot\frac 12 \alpha_2=2\left (\tan\frac 12 \beta_2+\tan\frac 12\gamma_1\right ).
\end{equation}

In the following case we do not give all small details of the calculations as those are very similar to the ones discussed above. We rather just point out the main conclusions of these calculations. 

Next, assume that for the side $x_1x_2$ it holds that $x_1\in K$ and $x_2\notin K$. Rotate the line $x_1x_2$ around $x_1$ by a sufficiently small positive angle $\varphi$ such that the intersection point $x_2'$ of the rotated line with the side $x_2x_3$ is still outside $K$. Let $P'$ be the polygon with vertices $x_1x_2'x_3\ldots x_n$. Clearly, $P'\subset P$. Let $b_2$ be the last intersection point of the side $x_1x_2$ with $\bd K$, as before. 

If the line $x_1x_2$ is not a supporting line of $K$ at $x_1$, then we obtain by hyperbolic trigonometry that 
\begin{multline*}
\dev (K,P')=\dev (K,P)+\varphi (2\sinh d(x_1,b_2)(\csc \beta_2-\cot\beta_2)\\-(\csc\alpha_2+\cot\alpha_2)\sinh d(x_1,x_2) )+O(\varphi^2) \text{ as }\varphi\to 0^+.
\end{multline*} 
Due to the optimality of $P$, it must hold that
\begin{equation}
\label{first}
\frac{\sinh d(x_1,x_2)}{\sinh d(x_1,b_2)}\cot\frac 12\alpha_2=2\tan \frac 12\beta_2.
\end{equation}
Note that $d(x_1,x_2)> d(x_1,b_2)$, and thus by the strictly monotonically increasing property of the $\sinh$ function it follows that the coefficient of $\cot (\alpha_2/2)$ in \eqref{first} is larger than $1$. 

If $x_1x_2$ is a supporting line of $K$, then, using \eqref{ferde-hossz}, we get that
\begin{multline*}
	\dev (K,P')=\dev (K,P)-\varphi(\csc\alpha_2+\cot\alpha_2)\sinh d(x_1,x_2)+O(\varphi^2) \text{ as }\varphi\to 0^+.
\end{multline*} 
As the coefficient of $\varphi$ is negative, this clearly contradicts the minimality of $P$, so $P$ cannot have such a side.

Now, the proof can be finished as in \cite[cf. (24)--(25) on p. 357]{Eggleston-1957}: For each $x_i\notin K$, the angle $\alpha_i$ appears in exactly two equations of type \eqref{second} or \eqref{first}, and $\beta_i$ and $\gamma_i$ in exactly one such equation.  Thus, by adding the two equations in which $\alpha_i$ appears, the coefficient of $\cot (\alpha_i/2)$ will be at least $2$. If $2+\varepsilon_i$ denotes the coefficient of $\cot (\alpha_i/2)$, then summing all equations of type \eqref{second} and \eqref{first} yields that
\begin{align*}
\sum (2+\varepsilon_i) \cot\frac 12 \alpha_i&=\sum 2\left (\tan\frac 12 \beta_i+\tan\frac 12\gamma_i\right )\\
&<\sum 2\tan\frac 12 (\beta_i+\gamma_i)\\
&\leq \sum 2 \tan \left (\frac{\pi}{2}-\alpha_i\right )\\
&=\sum 2\cot \frac 12\alpha_i, 
\end{align*}
which is clearly a contradiction as all $\varepsilon_i\geq 0$. 
This finishes the proof of Theorem~\ref{main}.

\section{Counterexample on the sphere}\label{sphere}
It is known that among spherical triangles contained in a (spherical) circle the inscribed regular triangle has the maximal perimeter, cf. L. Fejes T\'oth \cite{lagerungen}.
Thus, among triangles contained in the circle, the inscribed regular one has the minimum perimeter deviation from the circle. However, below we show an example of a triangle and circle, where the triangle is neither inscribed nor circumscribed, and approximates the circle better than either the inscribed or the circumscribed regular triangle.  

Let $K(r)$ be the spherical circle with centre $P$ and radius $r$. Consider the regular spherical triangle $T(d)=ABC\triangle$ with centre $P$ and circumradius $d$. Let $l=l(d)$ denote the side length of $T(d)$ and $m=m(d)$ its inradius, see Figure~\ref{Eggleston-ellenpelda}. 

\begin{figure}[h!]
	\begin{tikzpicture}[line cap=round,line join=round,>=triangle 45,x=1.0cm,y=1.0cm]
	\clip(5.16,-10.65) rectangle (14.32,-1.22);
	\draw [shift={(13.56,-8.39)}] plot[domain=2.09:3.14,variable=\t]({1*7.84*cos(\t r)+0*7.84*sin(\t r)},{0*7.84*cos(\t r)+1*7.84*sin(\t r)});
	\draw [shift={(5.72,-8.38)}] plot[domain=0:1.05,variable=\t]({1*7.84*cos(\t r)+0*7.84*sin(\t r)},{0*7.84*cos(\t r)+1*7.84*sin(\t r)});
	\draw [shift={(9.65,-1.6)}] plot[domain=4.19:5.23,variable=\t]({1*7.84*cos(\t r)+0*7.84*sin(\t r)},{0*7.84*cos(\t r)+1*7.84*sin(\t r)});
	\draw(9.63,-6.11) circle (3.85cm);
	\draw (9.63,-6.11)-- (13.42,-6.87);
	\draw (9.63,-6.11)-- (10.86,-2.46);
	\draw (9.63,-6.11)-- (9.65,-1.6);
	\draw (9.63,-6.11)-- (13.56,-8.39);
	\draw (10.34,-3.94) node[anchor=north west] {$r$};
	\draw (9.68,-3.32) node[anchor=north west] {$d$};
	\draw (9.63,-6.11)-- (12.56,-4.55);
	\draw (9.65,-8.71) node[anchor=north west] {$l$};
	\draw (6.81,-4.34) node[anchor=north west] {$l$};
	\draw (11.99,-4.78) node[anchor=north west] {$m$};
	\draw [shift={(9.63,-6.11)}] plot[domain=-0.2:0.49,variable=\t]({1*1.32*cos(\t r)+0*1.32*sin(\t r)},{0*1.32*cos(\t r)+1*1.32*sin(\t r)});
	\draw (10.45,-5.63) node[anchor=north west] {$\alpha$};
	\draw (12.65,-5.2) node[anchor=north west] {$s$};
	\begin{scriptsize}
	\fill [color=black] (5.72,-8.38) circle (1.5pt);
	\draw[color=black] (5.87,-8.14) node {$A$};
	\fill [color=black] (13.56,-8.39) circle (1.5pt);
	\draw[color=black] (13.7,-8.16) node {$B$};
	\fill [color=black] (9.65,-1.6) circle (1.5pt);
	\draw[color=black] (9.79,-1.36) node {$C$};
	\fill [color=black] (9.63,-6.11) circle (1.5pt);
	\draw[color=black] (9.77,-5.87) node {$P$};
	\fill [color=black] (13.42,-6.87) circle (1.5pt);
	\draw[color=black] (13.56,-6.63) node {$D$};
	\fill [color=black] (10.86,-2.46) circle (1.5pt);
	\fill [color=black] (12.56,-4.55) circle (1.5pt);
	\draw[color=black] (12.7,-4.31) node {$E$};
	\end{scriptsize}
	\end{tikzpicture}
\caption{\label{Eggleston-ellenpelda}}
\end{figure}

Then
$$\cos l=\cos^2 d+\sin^2 d\cos \frac{2\pi}{3},$$
and
$$\cos m=\frac{\cos d}{\cos \frac l2}.$$
Let $D$ be the intersection point of the side $BC$ with the circle $K(r)$ that is closer to $B$. Let $E$ be the intersection of the side $BC$ and the great circle through $P$ perpendicular to $BC$. Then $m$ is the distance of $P$ and $E$. Denote by $s=s(d)$ the length of the arc between $D$ and $E$. Furthermore, let $\alpha=\alpha(d)$ be the central angle $\angle EPD$.  Then
$$\cos s=\frac{\cos r}{\cos m},$$
and
$$\cos \alpha=\frac{\cos s-\cos m\cos r}{\sin m \sin r}.$$
Thus
$$f(r,d)=\dev_S (K(r), T(d))=6(2\alpha(d)\sin r-2s(d)+l(d)/2)-2\pi\sin r.$$
The graph of $f(\frac\pi2-0.1,d)$ over the interval $[\frac\pi 2-0.1, 1.52]$ is shown below, which clearly has its minimum inside the interval. We note that at the left endpoint of the interval the triangle is inscribed and at the right endpoint it is circumscribed. In fact, the circumscribed regular triangle approximates $K(r)$ better than the inscribed one, and the minimum occurs for a triangle that is neither inscribed nor circumscribed. Since all of these triangles are contained in the open hemisphere centred at $P$, they are convex in the spherical sense.

\begin{figure}[h!]
\includegraphics[scale=0.5]{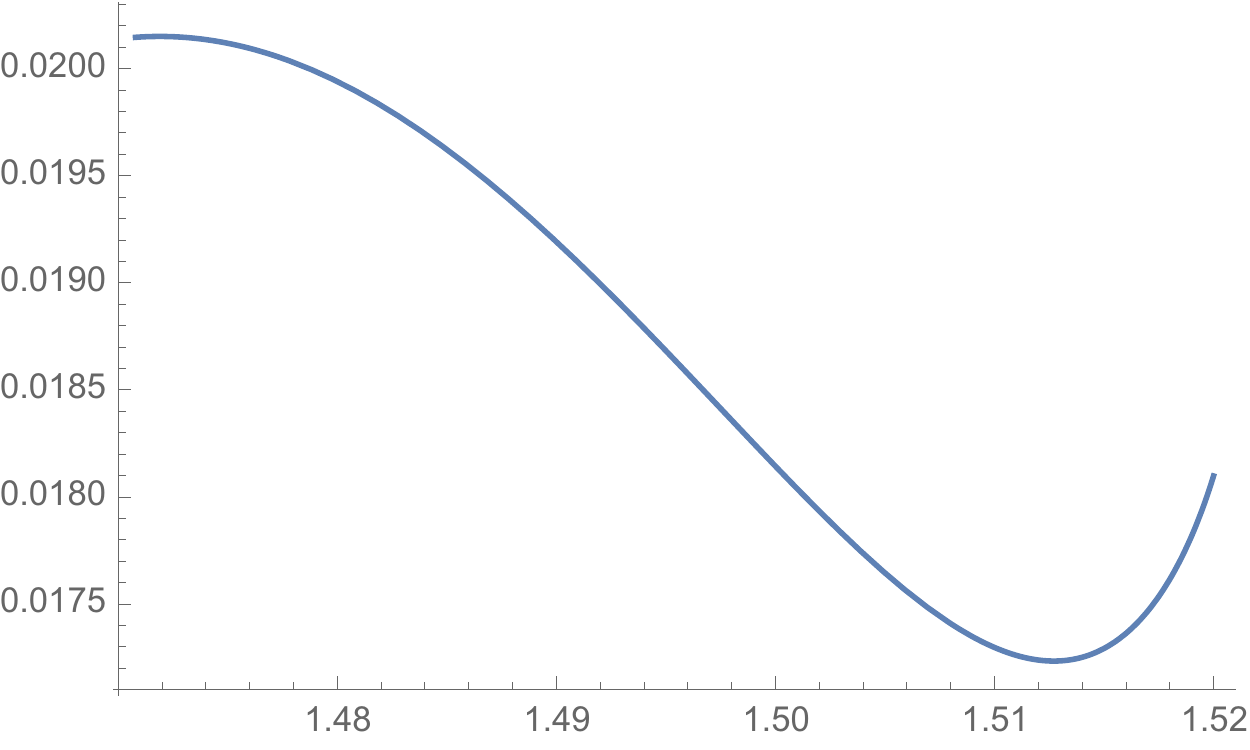}
\caption{The graph of $f(\frac\pi2-0.1,d)$ over the interval $[\frac\pi 2-0.1, 1.52]$}
\end{figure}

\section{Acknowledgements}
This research was partially supported by the National Research,
Development and Innovation Office of Hungary grant NKFIH K116451.

The author is grateful to Professor G\'abor Fejes T\'oth (Budapest, Hungary) for the enlightening discussions. The author also wishes to thank the MTA Alfr\'ed R\'enyi Mathematical Research Institute where part of this research was done.

\end{document}